\documentclass[12pt,a4paper,draft
]{article}
\usepackage{amsfonts,amssymb}
\usepackage{theorem}
{\theorembodyfont{\rm}
  \newtheorem{defi}{Definition}[section]
  \newtheorem{rem}[defi]{Remark}
  
  }
  \newtheorem{lem}[defi]{Lemma}
  \newtheorem{prop}[defi]{Proposition}
  \newtheorem{thm}[defi]{Theorem}
  \newtheorem{cor}[defi]{Corollary}
\newcommand{\Aa}{{\mathbb A}}

\newcommand{\CC}{{\mathbb C}}

\newcommand{\PP}{{\mathbb P}}

\newcommand{\RR}{{\mathbb R}}

\newcommand{\cB}{{\mathcal B}}

\newcommand{\cD}{{\mathcal D}}

\newcommand{\cR}{{\mathcal R}}
\newcommand{\cS}{{\mathcal S}}
\newcommand{\cT}{{\mathcal T}}
\newcommand{\cZ}{{\mathcal Z}}

\newcommand{\Aut}{{\mathrm{Aut}}}

\newcommand{\End}{{\mathrm{End}}}

\newcommand{\Hom}{{\mathrm{Hom}}}
\newcommand{\im}{{\mathrm{im}}}

\newcommand{\dis}%{{\mathrel{\vartriangle}}}
                 {{\mathrel{\scriptstyle{\triangle}}}}

\let\phi=\varphi

\let\theta=\vartheta
\newcommand{\SDelimArray}[4]{\hbox{\scriptsize\arraycolsep=.5\arraycolsep
  $\left#1\!\!\begin{array}{*{#3}{r}}#4\end{array}\!\!\right#2$}}

\newcommand{\SMat}{\SDelimArray()}

\newcommand{\qed}{\ $\square$}
\newcommand{\pf}{{\sl{Proof:}\ }}
\date{}
\begin{document}
\title{Affine Spaces within Projective Spaces}
\author{Andrea Blunck\thanks{Supported by a Lise Meitner
 Research Fellowship
of the Austrian Science Fund (FWF), project M529-MAT.} \and
{Hans Havlicek}}

\maketitle

%%%%%%%%%%%%%%%%%%%%%%%%%%%%%%%
\begin{abstract} \noindent
We endow the set of complements of a fixed subspace of a projective space
with the structure of an affine space, and
show that certain lines of such an affine space are affine reguli or cones
over affine reguli. Moreover, we apply our  concepts to the problem
of describing dual spreads.
We do not assume that the projective space is finite-dimensional
or pappian.

\noindent
{\em Mathematics Subject Classification} (1991):
51A30, 51A40, 51A45.\\
{\em Key Words:}  Affine space, regulus, dual spread.
\end{abstract}
%%%%%%%%%%%%%%%%%%%%%%%%%%%%%%%%%
\parskip .1cm
\parindent0cm
%%%%%%%%%%%%%%%%%%%%%%%%%%%%%%%%%%%%%%%%%%%%%%%%%%%%%%%%%%%%%%%%%%%
\section{Introduction}
%%%%%%%%%%%%%%%%%%%%%%%%%%%%%%%%%%%%%%%%%%%%%%%%%%%%%%%%%%%%%%%%%%%%
\thispagestyle{empty}

The aim of this paper is to equip the set of complements of a distinguished
subspace in an arbitrary projective space with the structure of an
affine space. In this introductory
section we first formulate the problem more
exactly and then describe certain special cases
that have already been treated
in the literature.

Let $K$ be a not necessarily commutative field, and let $V$ be a
left vector space over $K$ with arbitrary --- not necessarily finite
--- dimension. Moreover, fix a
subspace $W$ of $V$. We are interested
in the set
\begin{equation}\label{DefS}
\cS:=\{S\le W\mid V=W\oplus S\}
\end{equation} of all complements of $W$.

Note that often we shall take the projective point of view;
then we identify $\cS$ with the set of all those projective subspaces of the
projective space $\PP(K,V)$ that are complementary to the fixed
projective subspace induced by~$W$.

In the special
case that $W$ is a hyperplane, the set $\cS$ obviously is the point
set of an affine space, namely, of the
affine derivation of $\PP(K,V)$ w.r.t. this hyperplane.

This can be generalized
to arbitrary $W$:
In Section 2,  we are going
to introduce on $\cS$ the structure of a left vector space over $K$,
and hence the structure of an affine space. In general, this
affine space with point set $\cS$ will not be
determined uniquely.

In Sections 3 and 4 we study the lines of these affine spaces with point set
$\cS$ in terms of the projective space $\PP(K,V)$.
Finally, in Section 5 we use  our results in order
to describe dual spreads.

For the case that $K$ is commutative and $V$ is finite-dimensional over $K$,
the construction of the affine space on $\cS$ seems to be due to
J.G. Semple \cite{semp-31} (where $K=\RR,\CC$). Independently, in
\cite{metz-81} R. Metz shows that  then
$\cS$ can be identified with  the set of $(n-k)\times k$-matrices over~$K$
(where $n=\dim V$ and $k=\dim W$) and thus has the structure of a
vector space over~$K$ in a natural way.
The affine lines in $\cS$ are ``generalized reguli'',
i.e.,   reguli in  quotient spaces within  $\PP(K,V)$
(compare \cite{metz-81}).
In certain special cases this is  true in our more general situation
as well.

Metz's construction is also used by A. Herzer in \cite{herz-84}.
Moreover, in \cite{herz-84} and in \cite{semp-31}
a  representation of the affine space as a  ``stereographic
projection''  of the corresponding Grassmann variety is given.
For applications in differential geometry see
\cite{aki+g-97},   \cite{roux-91}, \cite{vinc-65}.

In
\cite{havl-87},  the  case of
lines skew to a fixed one  in  projective $3$-space
is generalized  to not necessarily pappian spaces via
 an explicit description
of the lines (cf. also \cite{havl-88b}).
In \cite{havl-87}
it is  also indicated how to generalize the used
approach to higher (but finite)
dimensions.
We will  not follow these ideas  here but  generalize
Metz's coordinatization
method and results on reguli obtained in \cite{blunck-99b}.

Finally we should mention  that one could also interpret
some results of this paper in terms of the extended concept of
chain geometry introduced in \cite{blu+h-99a}.

%%%%%%%%%%%%%%%%%%%%%%%%%%%%%%%%%%%%%%%%%%%%%%%%%%%%%%%%%%%%%%%%%%%
\section{Affine  spaces with point set $\cS$}
%%%%%%%%%%%%%%%%%%%%%%%%%%%%%%%%%%%%%%%%%%%%%%%%%%%%%%%%%%%%%%%%%%%%

Let  $W$ be a subspace of the  left vector
space $V$, and let $\cS$ be the set of complements of $W$  as introduced
in (\ref{DefS}). For our subsequent considerations we exclude the trivial
cases that
$W=\{0\}$ or $W=V$. In both cases  $\cS$ has only one element and thus
is a trivial affine space consisting of one single point.

Let $U\in\cS$ be a fixed complement of $W$.
We are going to coordinatize the
set  $\cS$ by the set $\Hom(U,W)$ of all $K$-linear mappings
$U\to W$, thus generalizing Metz's coordinatization by matrices.
Note that what follows is independent of the choice of $U$.

%We shall need the
%following  description of  the  ring $\End(V)$ of
% $K$-linear mappings $V\to V$:

Since $V=W\oplus U$,  each endomorphism $\phi\in\End(V)$
has the form
$$v=w+u\mapsto(w^\alpha+u^\gamma)+ (w^\beta+u^\delta),$$
($w\in W, u\in U$) with
linear mappings $\alpha:W\to W$, $\beta:W\to U$, $\gamma:U\to W$,
$\delta:U\to U$. So the endomorphism
ring $\End(V)$ is isomorphic to  the ring
$$\left\{\SMat2{\alpha&\beta\\\gamma&\delta}\mid\alpha\in\End(W),
\beta\in\Hom(W,U),\gamma\in\Hom(U,W),\delta\in\End(U)\right\},$$
equipped with the usual matrix
addition and multiplication (compare \cite{lang-95}, p.~643).
We shall frequently  identify $\End(V)$ with this matrix ring.

The stabilizer of our distinguished subspace $W$ in the group $\Aut(V)$
 is
\begin{equation}\label{stabil}
\left\{\SMat2{\alpha&0\\\gamma&\delta}\mid\alpha\in\Aut(W),
\gamma\in\Hom(U,W),\delta\in\Aut(U)\right\},
\end{equation}
where $0$ denotes the zero mapping $W\to U$.

Now consider an arbitrary complement of $W$, i.e., an element $S\in \cS$.
Because of $V=W\oplus U=W\oplus S$ there is a linear bijection
$\phi:V\to V$ fixing $W$ elementwise and mapping $U$ to $S$. This means
$\phi=\SMat2{1&0\\\gamma&\delta}$, where
$\gamma\in \Hom(U,W)$, $\delta\in \Aut(U)$, and $1=1_W$ is the identity mapping
on $W$. So
$$
S=U^\phi=\{u^\gamma+u^\delta\mid u\in U\}=:U^{(\gamma,\delta)}.$$
%\end{equation}
However, this description is not unique: We have
$U^{(\gamma,\delta)}=U^{(\gamma',\delta')}$ exactly if $(\gamma',\delta')=
(\rho\gamma,\rho\delta)$ for some $\rho\in\Aut(U)$.
In particular, since $\delta$ is invertible,
we may always assume $\delta=1=1_U$.

So we have
%\begin{equation}
$$\cS=\{U^{(\gamma,1)}\mid \gamma\in \Hom(U,W)\}$$
%\end{equation}
and we can  identify $\cS$ with $\Hom(U,W)$ via
\begin{equation}\label{ident}
U^{(\gamma,1)}\mapsto
\gamma.
\end{equation}

Note that instead of (\ref{ident}), we
 could also use the identification  $$\cS\ni U^{(\gamma,1)}=
 U^{(\rho\gamma,\rho)}\mapsto \rho\gamma\in\Hom(U,W),$$
for any  fixed $\rho\in \Aut(U)$.
One can easily check that this approach would yield
essentially  the same.

The abelian group $(\Hom(U,W),+)$ is a faithful
left module over the ring $\End(U)$ w.r.t.
composition of mappings.
Hence any embedding of some not necessarily commutative field $F$
into $\End(U)$ makes $\Hom(U,W)$ a left
vector space over $F$.
Analogously, $\Hom(U,W)$ is a faithful right module over $\End(W)$,
and any embedding $F\hookrightarrow\End(W)$ makes $\Hom(U,W)$ a right
vector space over $F$.

We will consider only the  case where $\Hom(U,W)$
becomes a left vector space over $K$. So we have to embed
the field
$K$ into the endomorphism ring $\End(U)$. In general, there are many
possibilities for such embeddings. We will
restrict ourselves to the  following  type:

\begin{lem}\label{embed}
\begin{enumerate}
\item
Let $(b_i)_{i\in I}$ be a basis of $U$. For each $k\in K$ we define
$\lambda_k\in\End(U)$ by $b_i\mapsto kb_i$ ($i\in I$). The mapping
 $\lambda:K\to
\End(U): k\mapsto \lambda_k$ is an injective homomorphism of rings, and
hence embeds the field $K$ into the ring $\End(U)$.
The associated left scalar multiplication  on $\Hom(U,W)$
is given by $k\cdot \gamma:=
\lambda_k \gamma\in\Hom(U,W)$ (for $k\in K$, $\gamma\in\Hom(U,W)$).
\item
Let $(b_i')_{i\in I}$ be another basis of $U$. Then the
 associated embedding $\lambda':K\to\End(U)$
 is given by $\lambda': k\mapsto
 \rho^{-1}\lambda_k\rho$, where $\rho\in\Aut(U)$ is the unique linear
 bijection with ${b_i}^\rho=b_i'$.
 \end{enumerate}
\end{lem}

Note that if $K$ is commutative, then each $\lambda_k$ is central in $\End(U)$,
and hence any two  embeddings of this type  coincide.
So in this case $\Hom(U,W)$ is a (left and right) vector space over $K$
in a canonical way. The same holds if we  consider $\Hom(U,W)$ as a
vector space over the center $Z$ of $K$.

The vector space structures on $\Hom(U,W)$ obtained in Lemma \ref{embed}
can be carried over to the set $\cS$ by using
the identification (\ref{ident}):

\begin{prop}\label{vspaces}
Let $(b_i)_{i\in I}$ be a basis of $U$, and let $\lambda:K\to\End(U)$
be the embedding associated to  $(b_i)_{i\in I}$ according to
Lemma \ref{embed}(a).

Then
$\cS$ is a  left vector space over $K$,  denoted by $(\cS,(b_i))$,
with addition
$$U^{(\gamma,1)}+U^{(\eta,1)}=U^{(\gamma+\eta,1)}$$
and scalar
multiplication
$$k\cdot U^{(\gamma,1)}=U^{(\lambda_k\gamma,1)}$$
(where $\gamma,\eta\in\Hom(U,W)$, $k\in K$).
\end{prop}

Each vector space
$(\cS,(b_i))$ gives rise to an  affine space with point set $\cS$.
This affine space  will be denoted by $\Aa(\cS,(b_i))$.

If the dimension of $V$ is finite, say $n$,  and $\dim W=k$, then
the affine space $\Aa(\cS,(b_i))$ is $k(n-k)$-dimensional.
If $\dim V=\infty$, then at least one of the subspaces $W$, $U$, and thus
also $\Aa(\cS,(b_i))$,
is infinite-dimensional (compare  Proposition \ref{PropA} below).

Now we consider the action of certain elements of
the stabilizer of $W$ in $\Aut(V)$ (i.e., the group~(\ref{stabil}))
on the affine spaces $\Aa(\cS,(b_i))$.

\begin{lem}\label{collin}
 \begin{enumerate}
 \item The group of all $\SMat2{1&0\\\eta&1}$ (with
 $\eta\in \Hom(U,W)$) induces on  $\cS$ the common translation group
 of all the affine spaces $\Aa(\cS,(b_i))$. This translation group is
 isomorphic to $(\Hom(U,W),+)$,  each translation has the form
 $U^{(\gamma,1)}\mapsto U^{(\gamma+\eta,1)}$.
 \item The group of all $\SMat2{\alpha&0\\0&1}$ (with $\alpha\in \Aut(W)$)
 induces   linear automorphisms of all vector spaces $(\cS,(b_i))$.
 So all these mappings are collineations of all $\Aa(\cS,(b_i))$, they
 have the form $U^{(\gamma,1)}\mapsto U^{(\gamma\alpha,1)}$.
 \item The group of all $\SMat2{1&0\\0&\rho}$ (with $\rho\in\Aut(U)$)
 induces linear bijections  $(\cS,(b_i))\to(\cS,({b_i}^\rho))$
 and hence isomorphisms between the associated affine spaces. In
 particular, this group
 acts transitively on the set of all  affine spaces  $\Aa(\cS,(b_i))$.
\end{enumerate}
\end{lem}
\pf
The mapping induced on $\cS$ by the matrix  $\SMat2{\alpha&0\\\eta&\rho}$
%(with $\alpha\in \Aut(W)$, $\eta\in\Hom(U,W)$, and $\rho\in \Aut(U)$)
is given by  $U^{(\gamma,1)} \mapsto U^{(\gamma\alpha+\eta,\rho)}
=U^{(\rho^{-1}(\gamma\alpha+\eta),1)}$. So assertions (a) and (b)
are obvious. Moreover, a closer look at Lemma  \ref{embed}(b) yields
assertion (c).
\qed

Note that the statements of
Lemma \ref{collin} (a) and (b) are valid for arbitrary affine spaces with
point set $\cS$ that are based upon the abelian group
$(\Hom(U,W),+)$, equipped with a scalar multiplication by embedding
$K$ into $\End(U)$: Of course all affine spaces of this type have
$(\Hom(U,W),+)$ as translation group. Moreover, $\gamma\mapsto \gamma\alpha$
($\alpha\in\Aut(W)$) is an endomorphism
of the left $\End(U)$-module $\Hom(U,W)$ and
hence induces a $K$-linear bijection w.r.t. all scalar multiplications under
consideration.

The permutations of $\cS$ that appear in Lemma \ref{collin}
are restrictions of
 projective
collineations of $\PP(K,V)$ fixing the subspace
 $W$. This will become important later when studying the
affine spaces $\Aa(\cS,(b_i))$ in terms of  $\PP(K,V)$.

In particular,
Lemma \ref{collin}(c) implies that any two of the  affine spaces
$\Aa(\cS,(b_i))$ are projectively equivalent.
Next we study under which conditions two of these affine spaces coincide.

\begin{lem}\label{coinci}
Let $(b_i)$ and $(b_i')$ be bases of $U$.
Then $\Aa(\cS,(b_i))=\Aa(\cS,(b_i'))$
 if, and only if, the unique element $\rho\in\Aut(U)$ with
 $b_i'={b_i}^\rho$ belongs to
$$N:=\{\nu\in\Aut(U)\mid \forall k\in K^*\
\exists l\in K^*:
\nu\lambda_k=\lambda_l\nu\},$$
which is the normalizer
of the image $(K^*)^\lambda$
of the multiplicative group $K^*$ under the embedding $\lambda$
in the group $\Aut(U)$.
\end{lem}
\pf
If the two affine spaces are affine lines, then they obviously coincide.
In this case, we have $\dim W=\dim U=1$, and hence in particular
$\Aut(U)=(K^*)^\lambda=N$.

Otherwise, the
 two affine spaces coincide exactly if the identity on $\cS$ is a
$K$-semilinear bijection $(\cS,(b_i))\to (\cS,(b_i'))$. Since both
vector spaces have the same additive group, this is equivalent to
the following:
\begin{equation}\label{semilin}
\exists \phi\in\Aut(K)\ \forall \gamma\in\Hom(U,W)\ \forall k\in K:
\lambda_k\gamma=\rho^{-1}\lambda_{k^\phi}\rho\gamma
\end{equation}
The $\End(U)$-module $\Hom(U,W)$ is  faithful, whence
(\ref{semilin}) is equivalent to
%\begin{equation}
$$\exists \phi\in\Aut(K)\ \forall k\in K:
\lambda_k=\rho^{-1}\lambda_{k^\phi}\rho.$$
%\end{equation}
This proves the assertion, because for each $\nu\in N$ the
mapping $k\mapsto l$ (where $\lambda_l=\nu\lambda_k\nu^{-1}$) is an automorphism
of $K$.
\qed

We want to describe the group $N$ in terms of the projective space.
This will yield
a projective characterization of the coincidence of two affine spaces of
type $\Aa(\cS,(b_i))$.

We need the following two  concepts:

Let $X$ be a  left
vector space over $K$, and let $(x_j)_{j\in J}$ be a basis
of $X$. The set of all linear combinations of the $x_j$'s with
coefficients in $Z$ is called the {\em $Z$-subspace of $X$ with respect
to $(x_j)$}. The associated set of points
of the projective space $\PP(K,X)$ is the  {\em (projective) $Z$-subspace of
$\PP(K,X)$ with respect to $(x_j)$}.

Now let $Y$ be another  left vector space over $K$, with basis
$(y_i)_{i\in I}$.
We call a linear mapping $\zeta:X\to Y$
 {\em central  w.r.t.  $(x_j)$
and $(y_i)$}, if it maps each $x_j$ into the  $Z$-subspace of $Y$
w.r.t. $(y_i)$. If in particular $X=Y$ and $(x_j)=(y_i)$, then
we say that $\zeta$ is {\em central w.r.t. $(x_j)$}.

Now we can formulate our next lemma.

\begin{lem}\label{descrN}
The group $N$ consists exactly of all products $\lambda_m\zeta$,
where $m\in K^*$ and $\zeta\in\Aut(U)$ is central w.r.t. $(b_i)$.

This means that $N$ consists of exactly those automorphisms of $U$ that
induce
projective collineations of $\PP(K,U)$ leaving the projective
$Z$-subspace w.r.t.
$(b_i)$ invariant.
\end{lem}
\pf
We only prove the first assertion, because the second one is just a
re-formulation.\\
 Obviously, all products $\lambda_m\zeta$ belong to $N$.
Conversely, consider $\nu\in N$. Fix one $k\in K^*$. Then there is an
$l\in k^*$ such that
\begin{equation}\label{Nbasis}
\forall i\in I: {b_i}^{\nu\lambda_k}=l{b_i}^\nu.
\end{equation}
Let $x$ be a coordinate (w.r.t. $(b_i)$) of any ${b_j}^\nu$. Then
(\ref{Nbasis}) implies $xk=lx$. Thus
if we fix one such non-zero coordinate, say $x_0$, then $x=x_0z_k$,
where $z_k$ centralizes $k$.
But this holds true for all $k\in K^*$, so that $x\in x_0Z$ for all
coordinates $x$ of all vectors ${b_j}^\nu$. Hence $\nu=\lambda_m\zeta$,
where $m=x_0$ and $\zeta$ is central with ${b_i}^\zeta=x_0^{-1}{b_i}^\nu$.
\qed

The following is a direct consequence of Lemmas \ref{coinci} and \ref{descrN}:

\begin{thm}\label{Zsub}
Let $(b_i)$ and $(b_i')$ be bases of $U$.
Then $\Aa(\cS,(b_i))=\Aa(\cS,(b_i'))$
if, and only if, the projective
$Z$-subspaces of $\PP(K,U)$ w.r.t. $(b_i)$ and
 $(b_i')$ coincide.
\end{thm}

This means that actually the affine space $\Aa(\cS,(b_i))$ does not depend
on the basis $(b_i)$ but only on the associated projective
$Z$-subspace of $\PP(K,U)$.

A special case is the following: If $\dim U=1$, then there is exactly one
projective $Z$-subspace in $\PP(K,U)$, namely, the point $U$ itself,
and hence there is a unique affine space of
type $\Aa(\cS,(b_i))$. Of course this is the affine derivation w.r.t.
the hyperplane $W$.

%%%%%%%%%%%%%%%%%%%%%%%%%%%%%%%%%%%%%%%%%%%%%%%%%%%%%%%%%%%%%%%%%%%
\section{The  symmetric case}
%%%%%%%%%%%%%%%%%%%%%%%%%%%%%%%%%%%%%%%%%%%%%%%%%%%%%%%%%%%%%%%%%%%%

As mentioned in the introduction,
we want to find out how the lines in the affine
spaces $\Aa(\cS,(b_i))$ look like in $\PP(K,V)$.
From now on we fix a basis $(b_i)$ and write $\Aa:=\Aa(\cS,(b_i))$
for short.
An arbitrary line of
$\Aa$ has the form
\begin{equation}\label{line}
\ell(\alpha,\beta):=\{U^{(\lambda_k\alpha+\beta,1)}\mid k\in K\}
\end{equation}
where $\alpha,\beta\in\Hom(U,W)$, $\alpha\ne 0$.

In this section we study  {\em symmetric case}, i.e., the case
that $W$ and its complement $U$ are isomorphic. In this special
situation we can use  results obtained in
\cite{blunck-99b}.
In analogy with
\cite{blunck-99b}, we assume w.l.o.g. that $V=U\times U$,
$W=U\times\{0\}$, and identify $U$  with $\{0\}\times U$.

In \cite{blunck-99b} the concept of a {\em regulus} in $\PP(K,V)$ is
introduced, generalizing the well known definition for the pappian case
(see, e.g., \cite{herz-95}) and also the definition for not necessarily pappian
$3$-space (cf.~\cite{segr-61}). It is shown in \cite{blunck-99b}
that each regulus is the image of the
{\em standard regulus}
%\begin{equation}\label{standreg}
$$\cR_0=\{W\}\cup\{U^{(\lambda_k,1)}\mid k\in K\}$$
%\end{equation}
under a projective collineation. In particular,  any two elements of a
regulus are complementary.

Obviously, the standard regulus equals the  line
$\ell(1,0)$ of $\Aa$, extended by $W$.
We call such a set $\cR\setminus\{W\}$, where $\cR$ is a regulus
containing $W$,
 an {\em affine regulus} (w.r.t. $W$). Each  affine regulus is
 a subset of $\cS$.

Moreover, we say that the line $\ell(\alpha,\beta)$ of $\Aa$ is
{\em regular}, if
$\alpha\in \End(U)$ is invertible. Since two elements $U^{(\gamma_1,1)}$
and $U^{(\gamma_2,1)}$ are complementary if, and only if, $\gamma_1-\gamma_2$
is invertible (compare~\cite{blunck-99a}), the regular lines
 are exactly the
lines joining two complementary
elements of $\cS$, and any two elements of such a line are complementary.

The following can be verified easily:

\begin{prop}\label{reglines}
A line of $\Aa$ is regular, exactly if it is an affine regulus
in $\PP(K,V)$.
In particular, the regular
line $\ell(\alpha,\beta)$ is the image of $\cR_0\setminus\{W\}$ under the
projective collineation induced by $\SMat2{\alpha&0\\\beta&1}$.
\end{prop}

Proposition \ref{reglines} implies that the affine space $\Aa$ possesses
a group of collineations that acts transitively on the set of regular lines,
because by Lemma \ref{collin} the matrix
$\SMat2{\alpha&0\\\beta&1}$ induces a collineation of $\Aa$.

Now we are looking for an explicit geometric description of the regular
lines as affine
reguli. First we need some more  information on
reguli.

Recall that a line
of $\PP(K,V)$ is called a {\em transversal}
of a regulus $\cR$, if it  meets each element of $\cR$ in exactly one point.
Through each point on a transversal of $\cR$ there
is a unique element of $\cR$.
The transversals of the standard regulus can be specified easily:

\begin{rem}\label{remTo}
The set $\cT_0$ of transversals of $\cR_0$
consists exactly of the lines $K(z,0)\oplus K(0,z)$, where
$z\ne 0$ belongs to the $Z$-subspace of $U$ w.r.t.
$(b_i)$.
\end{rem}

Recall, moreover, that for any three elements $U_1,U_2, U_3$ of a regulus
$\cR$ there is a perspectivity $\pi:U_1\to U_2$ with center $U_3$
(compare \cite{blunck-99b} for our general case).

\begin{lem}\label{transv}
Let $\cR$ be a regulus in $\PP(K,V)$, let $U_1, U_2,U_3\in\cR$ be pairwise
different, and let
$\cT$ be the set  of all transversals of $\cR$. Then the following holds:
\begin{enumerate}
\item The set $\cT\cap U_1:=\{T\cap U_1\mid T\in\cT\}$ is a projective
$Z$-subspace $\cZ(U_1)$
of $U_1$.
\item  If $\pi:U_1\to U_2$ is the  perspectivity with center $U_3$, then
 $\cZ(U_2)=\cZ(U_1)^\pi$. In particular,
$\cT=\{P\oplus P^\pi\mid P\in\cZ(U_1)\}$.
\item The regulus
$\cR$ is the set of all subspaces $X$  satisfying the following
conditions:
 \begin{enumerate}
 \item Each $T\in\cT$ meets $X$ in exactly one point, and $X$ is spanned
 by $\cT\cap X$.
 \item If $T_1\le T_2+T_3$ holds for $T_1,T_2,T_3\in\cT$, then the points
 $T_1\cap X$, $T_2\cap X$, $T_3\cap X$ are collinear.
 \end{enumerate}
\end{enumerate}
\end{lem}
\pf
The assertions (a) and (b) follow directly from
\ref{remTo}, because by \cite{blunck-99b}
the group $\Aut(V)$ acts
transitively on the set of triples $(\cR,U_1,U_2)$ where $\cR$ is a regulus
and
$U_1,U_2$ are different elements of $\cR$.\\
(c): Choose $T_1\in\cT$ and a point $P_1\in T_1$. Then there is a unique
$Y\in\cR$ through $P_1$. On the other hand, let $X$ be any subspace
through $P_1$ satisfying (i) and (ii). We have to show  that $X=Y$.\\
We first  observe that for each $T_2\in\cT\setminus\{T_1\}$ there is a
$T_3\in\cT\setminus\{T_1,T_2\}$ such that the point $T_3\cap  Y$ lies on the
line $L$ joining $T_1\cap Y$ and $T_2\cap Y$, since $\cT$ meets $Y$ in
a projective $Z$-subspace by (a).
Then the transversals $T_1,T_2,T_3$ all belong to
the $3$-space spanned by $L$ and $L^\pi$, where $\pi$ is the
perspectivity mapping $Y$ to another element of $\cR$,
according to (b). Moreover, also by (b),
any two transversals of $\cR$ are skew, and hence we have
$T_1\le T_2\oplus T_3$. Thus the points $T_i\cap X$ ($i=1,2,3$), that
exist by (i), are
collinear by (ii). Since there is a unique line through $P_1=T_1\cap X$
meeting
$T_2$ and $T_3$ (namely, the line $L$ joining $P_1$ and $T_2\cap Y=
(P_1\oplus T_3)\cap T_2$), we obtain $T_i\cap X=T_i\cap Y$. Applying this to
each $T_2\in\cT\setminus\{T_1\}$ yields  $X=Y$, because by (i)
$X$ is spanned by $\cT\cap X$.
\qed

Using Lemma \ref{transv}(c) one can  reconstruct the regulus $\cR$ from its set
of transversals.

\begin{cor}\label{reconst}
Every regulus is uniquely determined by the set of its
transversals. More exactly: If $\cR$ and $\cR'$ are reguli with the
same set  of transversals, then $\cR=\cR'$.
\end{cor}

Now consider a regulus $\cR$  containing $W$, with transversal set $\cT$.
We associate to it the
set
\begin{equation}\label{Tset}
W+\cT:=\{W+T\mid T\in\cT\}.
\end{equation}

\begin{cor}\label{idCor}
Let $U_1$ and $U_2$ be two different elements of $\cS$. Let $\cR$ and $\cR'$
be reguli containing $W$, $U_1$, $U_2$, and let $\cT$ and $\cT'$ be
their sets of transversals,
respectively. Then  $W+\cT=W+\cT'$  implies $\cR=\cR'$.
\end{cor}
\pf
Each transversal $T\in\cT$ is the line joining the points $(W+T)\cap U_1$
and $(W+T)\cap U_2$. So our assumptions imply $\cT=\cT'$ and hence
$\cR=\cR'$, by Corollary \ref{reconst}.
\qed

From now on, let $\cZ(U)$ be the projective
$Z$-subspace of $U$ w.r.t. $(b_i)$. Recall that by Theorem~\ref{Zsub} it is
$\cZ(U)$ rather than $(b_i)$ that determines the affine space $\Aa$. We set
%\begin{equation}\label{Zset}
$$W+\cZ(U):=\{W+P\mid P\in\cZ(U)\}.$$
%\end{equation}
The following is clear by  \ref{remTo}:
\begin{lem}\label{WplusTo}
Let $\cT_0$ be the set of transversals of  $\cR_0$.
Then $W+\cT_0=W+\cZ(U)$.
\end{lem}

Now we turn back to  the affine reguli that are lines of $\Aa$.
Note that by Proposition \ref{reglines}
such a line is necessarily regular.

\begin{thm}\label{trsvreg}
Let $\cR$ be a regulus containing $W$, and let $\cT$ be its set of
transversals. Then  $\cR\setminus\{W\}$ is a regular line of $\Aa$,
exactly if
\begin{equation}\label{WplusZ}
W+\cT=W+\cZ(U).
\end{equation}
\end{thm}
\pf
%``$\Rightarrow$'':
Let
$\cR\setminus\{W\}$ be the regular line $\ell(\alpha,\beta)$. Then
Proposition \ref{reglines} implies  $\cR={\cR_0}^\phi$, with $\phi$
induced by $\SMat2{\alpha&0\\\beta&1}$. Hence also  $\cT={\cT_0}^\phi$,
and we compute $W+\cT=W+{\cT_0}^\delta=(W+\cT_0)^\delta=W+\cT_0$.
Now Lemma \ref{WplusTo}
yields (\ref{WplusZ}).\\
%``$\Leftarrow$'':
Conversely, let $\cR$ be a regulus as in the assertion.
Choose complementary $U_1,U_2\in\cR\setminus\{W\}$. They are
joined by a regular line $\ell(\alpha,
\beta)$,
which %by the first direction
is  an affine  regulus $\cR'\setminus\{W\}$
with  transversal set
 satisfying
(\ref{WplusZ}). So Corollary \ref{idCor} yields $\cR=\cR'$ and
thus  $\cR\setminus\{W\}=\ell(\alpha,\beta)$.
\qed

\begin{cor}\label{detline}
Let $U_1, U_2\in\cS$ be complementary. Then the regular  line in $\Aa$
joining $U_1$ and $U_2$ is $\cR\setminus\{W\}$, where $\cR$
is the unique regulus containing $W$, $U_1$, $U_2$ that has a
transversal set
$\cT$ satisfying $W+\cT=W+\cZ(U)$.
\end{cor}

%%%%%%%%%%%%%%%%%%%%%%%%%%%%%%%%%%%%%%%%%%%%%%%%%%%%%%%%%%%%%%%%%%%
\section{The general case}
%%%%%%%%%%%%%%%%%%%%%%%%%%%%%%%%%%%%%%%%%%%%%%%%%%%%%%%%%%%%%%%%%%%%

Now we turn back to the general case that $V=W\oplus U$ for arbitrary
$W$ and~$U$.
We want to investigate the lines in $\Aa=\Aa(\cS,(b_i))$.
They all have the form $\ell(\alpha,\beta)$ as introduced in (\ref{line}).

It is sufficient
to consider the lines containing $U=U^{(0,1)}$, i.e., the lines of type
$\ell(\alpha,0)$, because the other lines are images of these under the
translation group introduced in Lemma \ref{collin}(a). More exactly,
the line
$\ell(\alpha,0)$ is mapped by $\SMat2{1&0\\\beta&1}$ to the
line $\ell(\alpha,\beta)$.
These  translations can also be considered
as collineations of the projective space $\PP(K,V)$ that fix all points
of  and all subspaces through the projective subspace induced by $W$.
On the other hand, each projective collineation fixing
  all points
of  and all subspaces through $W$ and mapping $U$ to $U^{(\beta,1)}$
is induced by a matrix $\SMat2{z&0\\\beta&1}$ with $z\in
Z^*$ and hence maps $\ell(\alpha, 0)$ to $\ell(\alpha,\beta)$.
Thus we have obtained:

\begin{lem}
The line $\ell(\alpha,\beta)$ of $\Aa$ is the image of $\ell(\alpha,0)$
under each projective collineation of $\PP(K,V)$ that
fixes all points
of  and all subspaces through $W$ and maps $U$ to $U^{(\beta,1)}$.
\end{lem}

We need the following slight generalization of Lemma \ref{collin}(b):

\begin{lem}\label{DefDach}
Let $V_1=W_1\oplus U$ and  $V_2=W_2\oplus U$ be  left vector spaces over $K$.
Let $\cS_j$ be the set of complements of $W_j$ in $V_j$,  and let
$\Aa_j=\Aa(\cS_j,(b_i))$ be the associated affine space ($j=1,2$).
If $\alpha:W_1\to W_2$ is  linear, then also
%\begin{equation}
$$\hat\alpha:V_1\to V_2: w_1+u\mapsto w_1^\alpha+u$$
%\end{equation}
is linear. The induced mapping  on  $\cS_1$, given by
\begin{equation}\label{Hat}
U^{(\gamma,1)}\mapsto U^{(\gamma\alpha,1)},
\end{equation}
and also denoted by $\hat\alpha$, is a
linear mapping $(\cS_1,(b_j))\to(\cS_2,(b_j))$ and hence
a  {\em homomorphism
of affine spaces}
$\Aa_1\to\Aa_2$, i.e., it maps lines to lines
or points and preserves parallelity.
\end{lem}

We apply this to $V_1=U\times U$ with $W_1=U\times \{0\}$, $U=\{0\}
\times U$, and $V_2=V=W\oplus U$, and thus obtain from
Proposition  \ref{reglines}:

\begin{prop}\label{project}
The line $\ell(\alpha,0)$  is the image of the standard affine
regulus  $\cR_0\setminus\{(U\times\{0\})\}$
in $\PP(K,U\times U)$ under the linear mapping $\hat\alpha:U^{(\gamma,1)}
\mapsto U^{(\gamma\alpha,1)}$.
\end{prop}

Since the linear mapping $\hat\alpha$ is defined on the whole vector space
$V_1=U\times U$, and thus on $\PP(K,U\times U)$, we can make a statement
on the image of the transversals of $\cR_0$ under $\hat\alpha$.

\begin{prop}\label{HatTrsv} Let $\ell(\alpha,0)$ be a line of $\Aa$ and let $\cR_0$
be the standard regulus in $\PP(K,U\times U)$. Then the  set $\cT_0$
of transversals of $\cR_0$ is mapped by $\hat\alpha$ onto a set
$\cT={\cT_0}^{\hat\alpha}$ of points and lines in $\PP(K,V)$, such that
the following holds:
\begin{enumerate}
\item If $T\in\cT$ is a point, then each element of $\ell(\alpha,0)$
contains $T$.
\item If $T\in\cT$ is a line, then $T$ is a {\em transversal}
of $\ell(\alpha,0)$,
i.e., $T$ meets $W$ and the mapping
$\ell(\alpha,0)\cup\{W\}\to T: X\mapsto X\cap T$ is a bijection.
\end{enumerate}
\end{prop}
\pf This is a direct consequence of the fact that $\hat\alpha$ is linear
and each restriction $\hat\alpha|_{U^{(\gamma,1)}}:U^{(\gamma,1)}
\to U^{(\gamma\alpha,1)}$ is a bijection.
\qed

One can also compute the set $\cT={\cT_0}^{\hat\alpha}$ explicitly: By
\ref{remTo}, each
transversal
$T\in\cT_0$ has the form $T=K(z,0)\oplus K(0,z)$ for some $z\ne0$ in the
$Z$-subspace of $U$ w.r.t. $(b_i)$. Its image is $T^{\hat\alpha}
=Kz^\alpha+Kz$, which is a point exactly if $z\in\ker(\alpha)$.

The  points of the line $\ell(\alpha,0)$  are
the subspaces
$U^{(\lambda_k\alpha,1)}=\{u^{\lambda_k\alpha}+u\mid u\in U\}$ with
$k\in K$. In particular, they all belong to the projective
space over $\im(\alpha)\oplus U\le V$. Using Lemma~\ref{DefDach}, we can
interpret this as follows:

\begin{rem}\label{Remark}
The line $\ell(\alpha,0)$ is a line in the affine space
$\Aa_0=\Aa(\cS_0,(b_i))$ of complements
of $\im(\alpha)$ in $\im(\alpha)\oplus U$. The affine space
$\Aa_0$  is embedded
into $\Aa$ via $\hat\iota$, where $\iota$ is the inclusion
mapping $\im(\alpha)\hookrightarrow W$. We will not
distinguish between $\ell(\alpha,0)$ as a line in $\Aa_0$ and
$\ell(\alpha,0)$ as a line in~$\Aa$.
\end{rem}

If $\alpha$ is injective, then the mappings
$\alpha:U\to \im(\alpha)$ and
$\hat\alpha:U\times U\to\im(\alpha)\oplus U$
are linear bijections.
In particular, $\hat\alpha$ induces an isomorphism of projective
spaces $\PP(K,U\times U)\to\PP(K,\im(\alpha)\oplus U)$,
and  $\ell(\alpha,0)$
is the image
of  the standard affine  regulus in $\PP(K,U\times U)$ under this
isomorphism.
 This, together with  Corollary \ref{detline}, yields the following:

\begin{prop}\label{invert}
If $\alpha\in\Hom(U,W)$ is injective, then $\ell(\alpha,0)$ is an
affine  regulus
in the projective subspace $\PP(K,\im(\alpha)\oplus U)$ of $\PP(K,V)$.\\
The associated regulus $\ell(\alpha,0)\cup \{\im(\alpha)\}$ is the
unique regulus in $\PP(K,\im(\alpha)\oplus U)$  containing
$\im(\alpha)$, $U$, and  $U^{(\alpha,1)}$  that has a transversal
set $\cT$ satisfying $\im(\alpha)+\cT=\im(\alpha)+\cZ(U)$.
\end{prop}

We need another type of homomorphism between affine spaces $\Aa(\cS,(b_i))$:

\begin{lem}\label{DefStern} Let $V_1=W\oplus U_1$ and
$V_2=W\oplus U_2$  be left vector spaces over $K$. Let
 $\cS_j$ be the set
of complements of $W$ in $V_j$, let $(b_i^{(j)})$ be a basis of $U_j$,
 and let  $\Aa_j=\Aa(\cS_j,(b_i^{(j)}))$ be the associated affine
 space ($j=1,2$).

If $\delta:U_1\to U_2$ is central w.r.t. $(b_j^{(1)})$
and $(b_j^{(2)})$ , then
%\begin{equation}
$$\delta^*: \cS_2\to\cS_1: U_2^{(\eta,1)}\mapsto U_1^{(\delta\eta,1)}$$
%\end{equation}
is a linear mapping $(\cS_2,(b_i^{(2)}))\to (\cS_1,(b_i^{(1)}))$
and hence a homomorphism of affine spaces $\Aa_2\to\Aa_1$.
\end{lem}
\pf
 Let $\lambda^{(j)}$ be the embedding of $K$ into $\End(U_j)$
w.r.t. $(b_i^{(j)})$. Using that $\delta$ is central, we obtain
$\lambda^{(2)}_k\delta=\delta\lambda^{(1)}_k$ for each $k\in K$. This
implies that $\delta^*$ is linear.
\qed

 If $\delta:U_1\to U_2$ and $\rho:U_2\to U_3$ are central linear mappings
w.r.t.  given bases, then the associated homomorphisms
of affine spaces obviously satisfy the condition
\begin{equation}\label{SternProd}
(\delta\rho)^*=\rho^*\delta^*.
\end{equation}

%Note that if $V$ is finite-dimensional, then one can obtain
%the linear mappings $\hat\alpha$
%(see (\ref{Hat})) and $\delta^*$ from the isomorphism of the left
%vector space $\Hom(U,W)$ with the tensor product ${}_K{U^*}_K\otimes{}_KW_Z$
%of bimodules (Zitat?). Here the dual vector space $U^*$ is considered
%as a left vector space over $K$ w.r.t. the embedding $\lambda:K\to\End(U)$
%associated to the basis $(b_i)$. Thus central linear mappings induce
%homomorphisms between the bimodules of type ${}_K{U^*}_K$.

We apply Lemma \ref{DefStern} and formula (\ref{SternProd})
to central subspaces of $U$. A subspace $U'\le U$
is called {\em central} (w.r.t. the fixed basis $(b_i)$),
if it possesses a basis of
$Z$-linear combinations of the $b_i$'s.
For the investigation of $\Aa$ it suffices (by Theorem \ref{Zsub})
to consider central subspaces with basis $(b_j)_{j\in J}$, where $J\subseteq
I$. Note that for each subspace of $U$ one can find a complement that is
central in $U$.

\begin{lem}\label{Stern}
Let $U'$ be a central subspace of $U$ with basis $(b_j)_{j\in J}$,
let $\cS'$ be the set of complements of $W$ in $W\oplus U'$,
and let $\Aa'=\Aa(\cS',(b_j))$ be the associated affine space. Then
the following
statements hold:
\begin{enumerate}
\item
The inclusion mapping $\iota:U'\hookrightarrow U$
is central (w.r.t. $(b_j)_{j\in J}$
and  $(b_i)_{i\in I}$).
The associated homomorphism $\iota^*:\Aa\to\Aa'$ is the ``intersection
mapping''
%\begin{equation}
$$U^{(\eta,1)}\mapsto U^{(\eta,1)}\cap (W\oplus U').$$
%\end{equation}

\item For any central complement $C$ of $U'$  the projection
$\pi:U\to U'$ with kernel $C$ is central. The associated homomorphism
$\pi^*:\Aa'\to \Aa$ is the ``join mapping''
%\begin{equation}
$$U'^{(\eta',1)}\mapsto U'^{(\eta',1)}\oplus C.$$
%\end{equation}
\item  The mapping $\pi^*\iota^*$ is the identity on $\Aa'$. So
$\pi^*$ is injective and $\iota^*$ is surjective.
\end{enumerate}
\end{lem}
\pf
(a): The definition of $\iota^*$ yields
$(U^{(\eta,1)})^{\iota^*}%=U'^{(\iota\eta,1)}
=\{(u'^\eta,u')\mid u'\in U'\}=U^{(\eta,1)}\cap (W\oplus U')$.\\
(b): We compute
$(U'^{(\eta',1)})^{\pi^*}%=\{(u'+c)^{\pi\eta'},u'+c)\mid u'\in U',c\in C\}
=\{(u'^{\eta'},u'+c)\mid u'\in U', c\in C\}=U'^{(\eta',1)}\oplus C$.\\
(c): This follows from (\ref{SternProd}), since $\iota\pi=1_{U'}$.
\qed

The following corollary  will be
important later:
\begin{cor}\label{affSub}
Let $C$ be a central subspace of $U$. Then the set
%\begin{equation}
$$\cS/C:=\{S\in\cS\mid C\le S\}$$
%\end{equation}
 is an affine subspace of $\Aa$.
\end{cor}
\pf
Let $U'\le U$ be any central complement of $C$,
let $\Aa'$ be
the affine space associated to $W\oplus U'$, and
let $\pi$ be the projection onto
$U'$ with kernel $C$. By Lemma \ref{Stern}, the linear
injection
$\pi^*$ embeds the affine space $\Aa'$ into $\Aa$, and the image
${\Aa'}^{\pi^*}$ equals $\cS/C$.
\qed

Lemma \ref{Stern}(a)
 enables us to investigate the lines of $\Aa$ by means of their
intersections with subspaces $W\oplus U'$, where $U'$ is central in $U$.
For every $\ell(\alpha,0)$ we choose an appropriate maximal central $U'$
and investigate
%\begin{equation}
$$\ell(\alpha,0)\cap (W\oplus U'):=\{X\cap (W\oplus U')
\mid X\in \ell(\alpha,0)\}.$$
%\end{equation}
By $\cZ(U')$ we denote the projective $Z$-subspace of $\PP(K,U')$ w.r.t.
the central basis of $U'$.

\begin{thm}\label{intersection}
Let $\ell(\alpha,0)$ be a line of $\Aa$, and let $U'\le U$ be a
central complement of $\ker(\alpha)$. Then the intersection
$$\ell(\alpha,0)
\cap (W\oplus U')=\ell(\alpha,0)\cap (\im(\alpha)\oplus U')$$
is an affine
regulus  in $\PP(K,\im(\alpha)\oplus U')\cong\PP(K,U'\times U')$.
The associated regulus $\cR$ is the unique one containing $\im(\alpha)$,
$U'$, $U'^{(\alpha,1)}$  that has a
transversal set $\cT$
satisfying $\im(\alpha)+\cT=\im(\alpha)+\cZ(U')$.
\end{thm}
\pf By Lemma \ref{Stern}(a) and Remark \ref{Remark}, the set $\ell(\alpha,0)
\cap (\im(\alpha)\oplus U')$ is a line in the affine space
$\Aa'$ associated to $\im(\alpha)\oplus U'$, namely $\ell(\iota\alpha,0)$,
where $\iota:U'\hookrightarrow U$ is
the inclusion. Since $U'$ is complementary to $\ker(\alpha)$,
we have that $\iota\alpha:U'\to \im(\alpha)$ is a bijection, whence
$\ell(\iota\alpha,0)$ is a regular line in $\Aa'$.
The rest follows from Corollary \ref{detline}.
\qed

Of course a similar statement holds for each central subspace $U'$ that
is skew to $\ker(\alpha)$. Then the intersection of $\ell(\alpha,0)$
with $U'^\alpha\oplus U'$ is an affine
regulus in $\PP(K,U'^\alpha\oplus U')$.
A special case is $U'=Kz$, with central $z\in U\setminus\ker(\alpha)$,
which leads to a transversal of $\ell(\alpha,0)$ as described in
Proposition \ref{HatTrsv}(b).

Now we consider  the kernel of $\alpha$.
Note that $\ker(\alpha)$ --- like any other subspace of $U$ ---
 contains a unique maximal central subspace.

\begin{thm}\label{KerMax}
Let $\ell(\alpha,0)$ be a line of $\Aa$, and let $M$ be the maximal central
subspace of $\ker(\alpha)$. Moreover, let $U'$ be a central complement of
$\ker(\alpha)$ and let $\cR$ be
the regulus in $\PP(K,\im(\alpha)\oplus U')$ appearing in
Theorem \ref{intersection}. Then the following
statements hold:
\begin{enumerate}
\item The set $\ell(\alpha,0)$ is entirely contained in
$$\cS/M=\{S\in\cS\mid M\le S\}.$$
\item The intersection $\ell(\alpha,0)\cap(\im(\alpha)\oplus U'\oplus M)$
 is the  {\em cone with  vertex $M$ over} the affine regulus $\cR\setminus
\{\im(\alpha)\}$, i.e., the set
$$\{X\oplus M\mid X\in\cR, X\ne \im(\alpha)\}.$$
\end{enumerate}
\end{thm}
\pf (a): By  Corollary
\ref{affSub} the set $\cS/M$ is an affine subspace of $\Aa$.
Since $M\le\ker(\alpha)=U\cap U^{(\alpha,1)}$, the line $\ell(\alpha,0)$
joining $U$ and $U^{(\alpha,1)}$ belongs to $\cS/M$.\\
(b): This follows directly from (a), Lemma \ref{Stern}(a), and
Theorem \ref{intersection}.
\qed

The bigger the maximal subspace $M$ of $\ker(\alpha)$ is, the better we know
$\ell(\alpha,0)$. In particular, if $\ker(\alpha)$ itself is central, we
obtain the following:

\begin{cor}\label{KerCent}
Let $\ell(\alpha,0)$ be a line of $\Aa$, and
let $\cR$ be
the regulus of Theorem \ref{intersection}. Assume  that $\ker(\alpha)$
is a central subspace of $U$. Then
$$\ell(\alpha,0)=\{X\oplus \ker(\alpha)\mid X\in\cR, X\ne\im(\alpha)\},$$
i.e.,  $\ell(\alpha,0)$ is the  cone with vertex $\ker(\alpha)$
over the affine regulus $\cR\setminus\{\im(\alpha)\}$.
\end{cor}

The statement of Corollary \ref{KerCent} means in other words that
``modulo $\ker(\alpha)$''
the line $\ell(\alpha,0)$ is an affine regulus.
This is what Metz in \cite{metz-81} called a generalized regulus. Since
Metz considers only commutative fields, in his case all lines of $\Aa$
are such generalized reguli.

The special case that  $\dim W=\dim U=2$ was treated in \cite{havl-87}.
Then each $\alpha\in\Hom(U,W)\setminus\{0\}$ that is not injective has
a one-dimensional kernel.
So there are only two types of non-regular
lines $\ell(\alpha,0)$, namely, those
where $\ker(\alpha)=Kz$ for some central
$z\ne0$, and the others.
In the first case we have a line pencil with carrier $Kz$ and one line
removed. This is the cone with vertex $Kz$ over an affine regulus, which
here is the set of points of an affine line.
%In the second case, the line $\ell(\alpha,0)$ is
%a degenerate conic of lines (see \cite{havl-87}).

\section{An application: Dual spreads}

In this section we aim at a description of dual spreads.

In a first step, where $V=W\oplus U$ is arbitrary, we show that
$(\cS,(b_i)_{i\in I})$ is isomorphic to the left vector space $W^I $ of all
$I$-families in $W$.

\begin{prop}\label{PropA}
The mapping
\begin{equation}\label{psi}
\psi:W^I\to\cS: (w_i)_{i\in I}\mapsto U^{(\gamma,1)},
\end{equation}
where $\gamma\in\Hom(U,W)$ is given by $b_i\mapsto w_i$, is a
$K$-linear bijection.
\end{prop}
\pf
Let $(w_i)^\psi=U^{(\gamma,1)}$ and $(w_i')^\psi=U^{(\gamma',1)}$.
Since  ${b_i}^{\gamma+\gamma'}=w_i+w_i'$ and ${b_i}^{\lambda_k\gamma}=kw_i$
hold  for all $i\in I$ and all $k\in K$, we have that $\psi$ is linear.
The bijectivity of $\psi$ is obvious.
\qed

This mapping $\psi$ has also a simple geometric interpretation:
For each $i\in I$ put
%\begin{equation}%\label{Ei}
$$E_i=W\oplus Kb_i.$$
%\end{equation}
Then $\PP(K,E_i)\setminus\PP(K,W)$ is an affine space isomorphic to
$\Aa(K,W)$, as follows from the isomorphism
\begin{equation}\label{psiI}
\psi_i:W\to \PP(K,E_i)\setminus\PP(K,W): w\mapsto K(w+b_i).
\end{equation}
Now $(w_i)^\psi$ is simply the element of $\cS$ that is spanned
by the points $K(w_i+b_i)$, where $i$ ranges in $I$.

Another tool for our investigation of dual spreads is the following:

\begin{lem}\label{lemA}
Let $(c_i)\in W^I $ and let $H\le W$ be a hyperplane. Then the $\psi$-image
of the affine subspace $(c_i)+H^I$ of $W^I$ equals
%\begin{equation}%\label{SofX}
$$\cS(X):=\{S\in\cS\mid S\le X\}$$
%\end{equation}
where $X:= (c_i)^\psi\oplus H$ is a hyperplane of $V$ with $W\not\le X$.\\
Conversely, each hyperplane $X\le V$ with $W\not\le X$ arises in this way
from exactly one affine subspace of $W^I$.
\end{lem}
\pf
Given $(c_i)$ and $H$ as above, then $X:=(c_i)^\psi\oplus H\le V$ is in fact
a hyperplane with $W\not\le X$, and $((c_i)+H^I)^\psi=\cS(X)$ follows from
(\ref{psiI}). \\
Conversely, if  $X\le V$ is a hyperplane with $W\not\le X$, then
$X$ contains an element of $\cS$, say $(c_i)^\psi$ with $(c_i)\in W^I$,
and $X\cap W=:H$ is a hyperplane of $W$. It is easily seen that
$((c_i)+H^I)^\psi=\cS(X)$.
\qed

In what follows we restrict ourselves to the symmetric case and
we adopt the settings of Section~3 ($V=U\times U$,
 $W=U\times \{0\}$, $U=\{0\}\times U$).
Recall that a line $\ell(\alpha,\beta)$ of $\Aa$ is regular if $\alpha\in
\Aut(U)$. Now a subspace of $\Aa$ will be called {\em singular},
if none of its
lines is regular.
We describe one family of maximal singular subspaces:

\begin{lem}\label{sing}
Let $X\le V$ be a hyperplane with $W\not\le X$. Then $\cS(X)=\{S\in\cS\mid
S\le X\}$ is a maximal singular subspace of $\Aa$.
\end{lem}
\pf
By Lemma \ref{lemA}, the set
$\cS(X)$ is a subspace of $\Aa$. Obviously, any two
elements of $\cS(X)$ are not complementary, whence
$\cS(X)$ is singular by Proposition \ref{reglines}. \\
Given $U_1\in\cS\setminus\cS(X)$ there exists a linear bijection
$\beta:W\to U_1$ such that $(W\cap X)^\beta=U_1\cap X$. Then
$Y:=\{w+w^\beta\mid w\in W\cap X\}$ is skew to $W\cap X$ and $U_1\cap X$.
There exists a complement $U_2\supseteq Y$ of $W\cap X$ relative to $X$.
By construction, $W$, $U_1$, and $U_2$ are pairwise complementary in $V$.
So the affine subspace spanned by $\cS(X)$ and $U_1$ cannot be singular.
\qed

A {\em dual spread} of $\PP(K,V)$ is a set of pairwise complementary
subspaces such that each hyperplane contains one of its elements.

\begin{prop}\label{propC}
A subset $\cB\subseteq\cS$ together with $W$ is a dual spread if,
and only if, the following conditions hold true:
\begin{description}
\item{\rm (DS1)} Distinct elements of $\cB$ are joined by a regular affine
line.
\item{\rm (DS2)} Each maximal singular subspace $\cS(X)$, where $X\le V$ is
a hyperplane with $W\not\le X$, contains an element of $\cB$,
\end{description}
\end{prop}
\pf
Let $\cB\subseteq \cS$. All elements of $\cB$ are complements of $W$,
and obviously every hyperplane of $V$ through $W$ contains $W\in\cB\cup
\{W\}$.\\
By Proposition \ref{reglines}, the elements of $\cB$ are pairwise
complementary exactly if {(DS1)} holds true. Moreover, {(DS2)}
is just a re-formulation of the condition that there is an element of $\cB$
in each hyperplane $X\le V$ with $W\not\le X$.
\qed

For another description of dual spreads containing $W$ we  generalize
the concept of   $*$-transversal mappings  described by N. Knarr
in \cite{knarr-95}, p.29. Moreover, we use Lemma \ref{lemA}, where
$W=U\times\{0\}$ can be replaced by the isomorphic vector space $U$.

A family $(\tau_i)_{i\in I}$ of mappings $\tau_i:D\to U$, where
$D\subseteq U$,
is called {\em $*$-transversal}, if the following conditions hold true:
\begin{description}
\item{\rm (T1*)} Given distinct $u,u'\in D$, then
$(u^{\tau_i}-u'^{\tau_i})_{i\in I}$ is a basis of $U$.
\item{\rm (T2*)} Given a family $(c_i)_{i\in I}\in U^I$ and a hyperplane
$H$ of $U$, there is a $u\in D$ such that $(u^{\tau_i})\in (c_i)+H^I$.
\end{description}

By {(T1*)} each $\tau_i$ is an injective mapping. So, after
an appropriate change of the domain $D$, one may always
assume that for one index $i_0\in I$ the mapping
$\tau_{i_0}$ is the canonical inclusion $D\hookrightarrow U$.
Therefore, in case $|I|=2$ it is enough to have one more mapping
$\tau_{i_1}:D\to U$. This is in fact the approach in \cite{knarr-95}.

%The subsequent theorem generalizes Knarr's result \cite{knarr-95}.

\begin{thm}
Let $(\tau_i)_{i\in I}$ be a $*$-transversal family of mappings $\tau_i:
D\to U$. Then
%\begin{equation}
$$\cD:=\{(u^{\tau_i})^\psi\mid u\in D\} \cup\{W\}=
\{\bigoplus\nolimits_{i\in I}
K(u^{\tau_i},b_i)\mid u\in D\} \cup\{W\}$$
%\end{equation}
is a dual spread of $\PP(K,V)$. Conversely, each dual spread of
$\PP(K,V)$ containing $W$ can be obtained in this way.
\end{thm}
\pf Let $(\tau_i)_{i\in I}$ be a $*$-transversal family. Choose distinct
elements
$(u^{\tau_i})^\psi= U^{(\gamma,1)}$ and $(u'^{\tau_i})^\psi=U^{(\gamma',1)}$,
where $u,u'\in D$ and $\gamma,\gamma'\in\End(U)$. Then
${b_i}^{\gamma-\gamma'}
=u^{\tau_i}-u'^{\tau_i}$ for all $i\in I$, so that $\gamma-\gamma'\in
\Aut(U)$ by  {(T1*)}. Hence $\cD\setminus\{W\}$ satisfies
{(DS1)}.
By {(T2*)} and Lemma~\ref{lemA}, $\cD\setminus\{W\}$ satisfies
{(DS2)}. So Proposition \ref{propC} shows that $\cD$ is a dual spread.\\
On the other hand, let $\cD$ be a dual spread containing $W$.
If  $(s_i),(s_i')\in (\cD\setminus\{W\})^{\psi^{-1}}$
coincide in one entry, say $s_j=s_j'$,
then $K(s_j,b_j)$ is a common point of their $\psi$-images by (\ref{psiI}).
Hence in this case $(s_i)=(s_i')$. \\
We fix one index $i_0\in I$ and put
$D:=\{s_{i_0}\mid (s_i)\in(\cD\setminus\{W\})^{\psi^{-1}}\}\subseteq U$.
By the above,  for each $u\in D$ there is a unique $(s_i)\in
(\cD\setminus\{W\})^{\psi^{-1}}$ with $u=s_{i_0}$. Moreover,
$\tau_i:D\to U: u\mapsto s_i$ is well defined for all $i\in I$. By
reversing the arguments of the first direction,
it is easily seen that the family
$(\tau_i)_{i\in I}$ is
$*$-transversal.
\qed
%%%%%%%%%%%%%%%%%%%%%%%%%%%%%%%%%%%%%%%%%%%%%%%%%%%%%%%%%%%%%%%%%%%%%%%%%%
{\small
%\bibliographystyle{plain}
%\bibliography{ketten}

\bigskip

Institut f\"ur Geometrie\\
Technische Universit\"at\\
Wiedner Hauptstra{\ss}e 8--10\\
A--1040 Wien,
Austria}
\end{document}